\documentclass{ws-ijgmmp}

\usepackage{amsfonts,amssymb}

\renewcommand{\theequation}{\thesection.\arabic{equation}}
\begin{document}

\markboth{ALEXANDER I. NESTEROV } {SMOOTH LOOPS AND FIBER BUNDLES:
THEORY OF PRINCIPAL Q-BUNDLES}

%
\catchline{}{}{}{}{}
%
\title{SMOOTH LOOPS AND FIBER BUNDLES: THEORY OF PRINCIPAL Q-BUNDLES}

\author{ALEXANDER I. NESTEROV }

\address{ Departamento de F{\'\i}sica, CUCEI, Universidad de
Guadalajara, Av. Revoluci\'on 1500\\ Guadalajara, CP 44420, Jalisco,
M\'exico\\
\email{nesterov@cencar.udg.mx } }

\maketitle

\begin{history}
\received{(Day Month Year)} \revised{(Day Month Year)}
\end{history}
\begin{abstract}
A nonassociative generalization of the principal fiber bundles
with a smooth loop mapping on the fiber is presented. Our approach
allows us to construct a new kind of gauge theories that involve
higher ''nonassociative'' symmetries.
\end{abstract}




\section{Introduction}

During the last few decades, nonassociative structures have been
employed in various fields of modern physics. Among others, one
may mention the rise of nonassociative objects such as
3-cocycles, which are linked with violations of the Jacobi identity in
anomalous quantum field theory, and quantum mechanics with the Dirac
monopole, the appearance of Lie groupoids and algebroids in the context
of Yang-Mills theories, and the application of nonassociative algebras
to gauge theories on commutative but nonassociative fuzzy spaces
\cite{Jac,G2,Gr,Gr1,G1,N1a,N3,NF,N1,NI,Strobl,Ram,MR}.

Nonassociative algebraic structures such as quasigroups and loops have
considerable potential interests for mathematical physics, especially in view
of the appearance of nonassociative algebras, such as the Mal'cev algebra
\cite{Mal}, related to the problem of the chiral gauge anomalies, the emergence of a
nonassociative electric field algebra in a two-dimensional gauge theory, and so
on \cite{C,Jo,NSem,Sem,Jac}. Quasigroups and loops have recently been
employed in general relativity and for the description of the Thomas precession,
coherent states, geometric phases, and nonassociative gauge theories
\cite{NI,N2,N4,N5,NS3,NS3a,NS2,NS1,U1,U2,K}.

In this paper, we give a detailed account of nonassociative fiber bundles
leading to ``nonassociative'' gauge theories \cite{N1,NI,N2}. In general terms,
a consequence of nonassociativity is that the structure constants of the gauge
algebra have to be changed by structure functions.

The paper is organized as follows. In Section 2, we outline the basic
constructions of the smooth loops theory. In Section 3,
the theory of nonassociative principal fiber bundles ($Q$-bundles)
is introduced.

\section{Smooth Loops}

\subsection{Basic notations }

Here, we outline the main results on the algebraic theory of quasigroups
and loops. The details may be found in \cite{Bel,Bruck,Pfl,Ch,Sab1}.

Let $\langle Q,{\mathbf \cdot}\rangle $ be a groupoid, i.e. a set
with a binary operation $(a,b) \rightarrow{} {a \mathbf \cdot b}$. A
groupoid $\langle Q,{\mathbf \cdot}\rangle $ is called a {\it
quasigroup} if each of the equations $a{\mathbf \cdot}
x=b,~y{\mathbf \cdot} a=b$ has unique solutions:  $x=a\backslash b$,
$y=b/a$. A {\it loop} is a quasigroup with a two-sided identity
$a{\mathbf \cdot} e= e{\mathbf \cdot} a=a, \forall a \in Q$. A loop
$\langle Q,{\mathbf \cdot},e\rangle $ with a smooth operation
$\phi(a,b):=a{\mathbf \cdot} b$ is called the {\it smooth loop}. We
define
\begin{eqnarray}
&L_a b=R_b a=a{\mathbf \cdot} b,\quad
l_{(a,b)}=L^{-1}_{a{\mathbf \cdot} b}\circ L_a\circ L_b, \nonumber\\
& \hat l_{(a,b)}= L_a\circ L_b \circ L^{-1}_{a{\mathbf \cdot} b} ,
\quad r_{(b,c)} = R^{-1}_{b{\mathbf \cdot}c}\circ R_c \circ R_b,
\label{Ll}
\end{eqnarray}
where $L_a$ is a {\it  left translation}, $R_b$ is a {\it right
translation}, $l_{(a,b)}$ is a {\it left associator}, $\hat
l_{(a,b)}$ is an {\it adjoint associator} and $r_{(b,c)}$ is a {\it right
associator}.

Let $T_e(Q)$ be the tangent space of $Q$ at the neutral element $e$. Then
for each $X_e \in T_e(Q)$, we construct a smooth vector field on $Q$
\[
X_b = L_{b\mathbf  \ast } X_e, \quad b\in Q, \quad X_e \in T_e(Q),
\; X_b \in T_b(Q)
\]
where $L_{b\mathbf \ast }: T_e(Q)\rightarrow{} T_b(Q)$ denotes the
differential of the left translation. Notice that $X_b$ depends smoothly on both variables $b \in Q$ and $X_e \in T_e(Q)$, and satisfies
\begin{eqnarray}
L_{a\mathbf  \ast }X_b = \hat l_{(a,b){\mathbf \ast} } X_{a{\mathbf
\cdot} b}.
\end{eqnarray}

\begin{definition}
A vector field $X$ on $Q$ that satisfies the
relation $L_{a\mathbf  \ast }X_b =\hat l_{(a,b)\mathbf  \ast }
X_{a{\mathbf \cdot} b}$ for any $a,b \in Q$ is called the {\it left
fundamental} or {\it  left quasi-invariant} vector
field.
\end{definition}

Let ${V}$ denote the set of left fundamental vector
fields on $Q$. It becomes a vector space under the operations
\begin{eqnarray}
&&(X + Y)_a = X_a + Y_a, \quad X,Y \in {V},\; a\in Q\nonumber,\\
&&(\alpha X)_a = \alpha X_a , \quad \alpha \in {\mathbb R}.
\nonumber
\end{eqnarray}
The following lemma is evident.

{\begin{lemma} The vector space of left fundamental vector fields is
isomorphic to the tangent space $T_e(Q)$ at the neutral element. The
isomorphism is defined by the map $L^{-1}_{a\mathbf  \ast}$:
$X_a\rightarrow{} X_e, \; X_e\in T_e(Q), \; X_a\in
T_a(Q)$.
\end{lemma}}

The vector space $V$ can be equipped with a Lie commutator. This leads us
to the notion of {\it quasialgebra}.

{\begin{definition} We define a {\it quasialgebra} $\mathfrak q$ on
$Q$ as the vector space of left fundamental vector fields, under the
Lie commutator operation.\end{definition}}

Let $\Gamma_i =R^j_i\partial/\partial a^j, \; i=1,2, \dots, r$ be a basis of the space of  left fundamental vector fields. Then, we have
\begin{eqnarray}
[\Gamma_i, \Gamma_j] = C^p_{ij}(a) \Gamma_p
\end{eqnarray}
where $C^p_{ij}(a)$ are the {\em structure functions} satisfying the
modified Jacobi identity
\begin{eqnarray}
C^p_{ij,n}R^n_k + C^p_{jk,n}R^n_i + C^p_{ki,n}R^n_j
+ C^l_{ij}C^p_{kl} +  C^l_{jk}C^p_{il} + C^l_{ki}C^p_{jl}=0.
\label{eq2}
\end{eqnarray}

In view of the noncommutativity of the right and left translations,
there exists a problem in the definition of the `adjoint' map of
$Q$ on itself. Indeed, in the case of the group $G$, the adjoint
map is given by
\[
Ad_g g' = g g' g^{-1} \equiv L_g \circ R^{-1}_g g' \equiv
 R^{-1}_g \circ L_g g'.
\]

For quasigroups, the left and right translations do not commute, $L_a
\circ R_b  \neq R_b \circ L_a $, and it follows that the definitions
$Ad_g g' = L_g \circ R^{-1}_g g'$ and $\tilde Ad_g g' =R^{-1}_g
\circ L_g g'$ are not equivalent. We introduce a generalized
adjoint map of $Q$ on itself in the following way.
{\begin{definition} The map
\begin{eqnarray}
Ad_b(a)= L^{-1}_a\circ R^{-1}_b\circ L_{a{\mathbf \cdot}b}: \;
Q\rightarrow{} Q \label{A_1}
\end{eqnarray}
is called the {\it Ad-map.} \end{definition}}

{\begin{remark} \rm The Ad-map $Ad_b(a)$ leaves invariant the
neutral element $e$ and generates the map $T_e(Q) \rightarrow{}
T_e(Q)$ as follows:
\begin{eqnarray}
{\sf Ad}_b(a):=\bigl(Ad_b(a)\bigr)_{\mathbf \ast}= L^{-1}_{a\ast}
R^{-1}_{b\ast} L_{a{\mathbf \cdot}b \ast}.
\end{eqnarray}
It is easy to check that $Ad_b(e)= Ad_b:=R^{-1}_b\circ
L_b$.\end{remark}}

{\begin{definition} The vector-valued 1-form
$\mbox{\boldmath$\omega$}$ defined through the relation
\begin{eqnarray}
\mbox{\boldmath$\omega$}(V_a)=V_e \quad {\rm where}  \quad V_e
=L^{-1}_{a^{\mathbf \ast} } V_a, \quad V_e \in T_e(Q), \; V_a \in
T_a(Q)
\end{eqnarray}
where $V$ is a left fundamental vector field, is called the {\em
canonical} $ \sf Ad$-form. \end{definition}}

{\begin{theorem} The canonical $\sf Ad$-form ${\omega}$ is a left
fundamental form and is transformed under left (right) translations as follows:
\begin{eqnarray}
&(L^{\mathbf \ast}_b\omega)V_a=l_{(b,a)\mathbf \ast}
\omega (V_a) \label{AL} \\
&(R^{\mathbf \ast} _b\omega)V_a={\sf Ad}^{-1}_b(a) \omega(V_a).
\label{AR}
\end{eqnarray}
\label{Th1}\end{theorem}}

{\begin{proof} For the left translations, we obtain
\[ (L^{\mathbf \ast} _b\omega)V_a=\omega(L_{b\mathbf  \ast
}V_a)=L^{-1}_{b{\mathbf \cdot} a\ast}L_{b\mathbf  \ast } L_{a\mathbf
\ast } V_e = l_{(b,a)\mathbf  \ast }\omega (V_a),
\]
thus (\ref{AL}) is true.

Similarly, for the right translations, we find
\begin{eqnarray}
(R^{\mathbf \ast} _b\omega)V_a=&\omega(R_{b\mathbf  \ast }V_a)
=L^{-1}_{a{\mathbf \cdot} b\mathbf  \ast }R_{b\mathbf  \ast }
L_{a\mathbf
\ast }L^{-1}_{a\mathbf  \ast } V_a \nonumber  \\
= &L^{-1}_{a{\mathbf \cdot} b\mathbf  \ast } R_{b\mathbf  \ast }
L_{a\mathbf  \ast }\omega (V_a) ={\sf Ad}^{-1}_b(a) \omega(V_a)
\end{eqnarray}
and the proof is completed. \end{proof}}

\subsubsection{Examples}

The examples given below play an important role in the description of generalized
coherent states, Thomas precession and nonassociative geometry \cite{N3,NS3,NS3a,NS2,NS1,U1,U2}.

\begin{example}
Loop $\mathbf \langle {\mathbb R}/{\mathbb
Z},\mathbf \ast\rangle $. Let $\ast$ be the operation defined on
${\mathbb R}/{\mathbb Z}$ as follows:
\begin{eqnarray}
&&x\ast y = x+y + f(x) +f(y) - f(x+y)
\end{eqnarray}
where $f(x)= (1-\cos 2\pi x)/4 $. Then $\langle {\mathbb
R}/{\mathbb Z}, \ast\rangle $ is the analytical loop \cite{H1,H2}.
\end{example}

\begin{example}
Loop $\rm Q\mathbb C \label{ejemplo}$. Let
${\mathbb C}$ be a complex plane and $\zeta,\eta\in{\mathbb C}$.
Then, the loop $\rm Q\mathbb C $ is obtained by introducing on
${\mathbb C}$ the nonassociative multiplication defined as follows
\begin{eqnarray}
R_\eta^\zeta=L_\zeta \eta=\zeta\mathbf
\cdot\eta=\frac{\zeta+\eta}{1-\overline\zeta\eta}, \quad
\zeta,\eta\in{\mathbb C}, \label{eq5}
\end{eqnarray}
where a bar denotes complex conjugation. The inverse operation is is found to
be
\begin{eqnarray}
L^{-1}_\zeta \eta=\frac{\eta-\zeta}{1+\overline\zeta\eta}, \quad
\zeta,\eta\in{\rm Q\mathbb C}. \label{eq6}
\end{eqnarray}
The associator $l_{(\zeta,\eta)}$ is given by
\begin{eqnarray}
l_{(\zeta,\eta)}\xi=\frac{1-\zeta\overline\eta}{1-\eta\overline\zeta}\xi
\end{eqnarray}
and can also be written as $l_{(\zeta,\eta)}= \exp({\rm i}\alpha),$ where $
\alpha =2\arg(1-\zeta\overline\eta)$.

The basis of left quasi-invariant vectors and the dual basis of 1-forms
are found to be
\begin{eqnarray}
\Gamma_1 = (1 + |\eta|^2)\partial_\eta,
\quad \Gamma_2 = (1 + |\eta|^2)\partial_{\bar\eta}, \\ \theta^1 =
\frac{d\eta}{1 + |\eta|^2}, \quad \theta^2 = \frac{d\bar\eta}{1 +
|\eta|^2}.
\end{eqnarray}
The computation of the commutator yields
\begin{eqnarray}
[\Gamma_1,\Gamma_2] = \bar\eta\Gamma_2 - \eta\Gamma_1,
\end{eqnarray}
And, for the structure functions, we find $C^1_{12}= -\eta, C^2_{12}=
\bar\eta$.

The loop thus obtained is related to the two-sphere $S^2$, which admits a natural quasigroup
structure. Namely, $S^2$ is the local two-parametric Bol loop $QS^2$
\cite{N3,N1,NS2,NS1}. The isomorphism between points of the sphere
and the complex plane ${\mathbb C}$ is established by the
stereographic projection from the south pole of the sphere:
$\zeta=e^{\rm i\varphi}\tan(\theta/2)$.
\end{example}

\begin{example}
Loop $QSU(2)$. Let us consider a set of the unitary
matrices $QSU(2)\subset$ $SU(2)$, where an arbitrary element of
$QSU(2)$ has the form
\begin{displaymath}
U=
\left ( \begin{array}{cc}
\alpha & \beta    \\
- \bar\beta & \bar\alpha
\end{array}
\right ), \quad |\alpha|^2 +  |\beta|^2 =1, \quad \alpha =\bar\alpha.
\end{displaymath}

It is convenient to set
\[
\alpha = \frac{1}{\sqrt{1+ |\eta|^2}}, \quad
\beta  = \frac{\eta}{\sqrt{1+ |\eta|^2}},
\]
and write $U$ as
\begin{eqnarray}
U_{\eta} = \frac{1}{\sqrt{1+ |\eta|^2}}
\left ( \begin{array}{cc}
1 & \eta    \\
- \bar\eta & 1
\end{array}
\right ).
\label{QS}
\end{eqnarray}

For arbitrary elements of $QSU(2)$, we define the nonassociative binary operation
as follows:
\begin{eqnarray}
U_{\eta}\ast U_{\zeta} = U_{\eta} U_{\zeta} \Lambda (\eta, \zeta)
\label{U},
\end{eqnarray}
where, on the right-hand side, the matrices are multiplied in the usual way, and
\begin{displaymath}
\Lambda (\eta, \zeta) =
\left ( \begin{array}{cc}
e^{i\varphi} & 0    \\
0           & e^{-i\varphi}
\end{array}
\right ),
\quad \varphi = 2\arg(1 - \bar\eta \zeta).
\end{displaymath}
The set of matrices (\ref{QS}) with the operation $\ast$ is called the loop
$QSU(2)$.

The loop $QSU(2)$ forms the so-called nonassociative representation of
$QS^2$. Indeed, writing  Eq. (\ref{U}) as
\begin{eqnarray}
U_{\eta}\ast U_{\zeta} = \frac{1}{\sqrt{1+ |\eta{\mathbf
\cdot}\zeta|^2}} \left ( \begin{array}{cc}
1 & \eta {\mathbf \cdot}\zeta   \\
- \overline{\eta {\mathbf \cdot}\zeta}& 1
\end{array}
\right ),
\label{QS1}
\end{eqnarray}
we find (see Ex. 9)
\begin{eqnarray}
L_{\eta}\zeta \equiv\eta{\mathbf
\cdot}\zeta=\frac{\eta+\zeta}{1-\overline\eta\zeta}, \quad
l_{(\eta,\zeta)}\xi =\frac{1- \eta\overline\zeta}
{1-\overline\eta\zeta}\xi, \quad \zeta,\eta,\xi\in{\mathbb C}.
\end{eqnarray}
\end{example}
\begin{example}
Loop $QH^2$. This loop is associated with the group
$SU(1,1)$ and its action on the two-dimensional unit hyperboloid
$H^2$. Let $D\subset\mathbb C$ be the open unit disk,
$D=\{\zeta\in\mathbb C:\left|\zeta\right|<1 \}$. We define the
nonasociative binary operation $\ast$ as
\begin{eqnarray}
\zeta\ast\eta=\frac{\zeta+\eta}{1+\bar\zeta\eta}, \quad\zeta,\eta\in D.
\end{eqnarray}
The associator $l_{(\zeta,\eta)}$ on $QH^2$ is determined by
\begin{eqnarray}
l_{(\zeta,\eta)}\xi=\frac{1+\zeta\overline\eta}{1+\eta\overline\zeta}\xi,
\end{eqnarray}
and can also be written as $l_{(\zeta,\eta)}=\exp({\rm i}  \alpha), \;
\alpha =2\arg(1+\zeta\overline\eta)$. Inside the disk $D$, the set of
complex numbers with the operation $\ast$ forms a two-sided loop $QH^2$, which
is isomorphic to the geodesic loop of two-dimensional Lobachevsky space
realized as the upper part of the two-sheeted unit hyperboloid. The
isomorphism is established by $\zeta=e^{\rm i\varphi}\tanh(\theta/2)$,
where $(\theta,\varphi)$ are inner coordinates on $H^2$.

The basis of left quasi-invariant vectors and the dual basis of 1-forms are
found to be
\begin{eqnarray}
\Gamma_1 = (1 - |\eta|^2)\partial_\eta, \quad
\Gamma_2 = (1 - |\eta|^2)\partial_{\bar\eta}, \\
\theta^1 = \frac{d\eta}{1 - |\eta|^2}, \quad
\theta^2 = \frac{d\bar\eta}{1 - |\eta|^2},
\end{eqnarray}
and the computation of the commutator gives
\begin{eqnarray}
[\Gamma_1,\Gamma_2] = \eta\Gamma_1 - \bar \eta\Gamma_2.
\end{eqnarray}
\end{example}

\begin{example} Loop $\rm Q{\sf H}_{\mathbb R}$. Let us consider
the quaternionic algebra over the complex field ${\mathbb C}(1,\rm i)$
\[ {\sf H}_{\mathcal C} = \{ \alpha + \beta i+ \gamma j + \delta k\;
\mid \; \alpha , \beta,\gamma,\delta \in
 {\mathbb C} \} \] with the multiplication operation defined by the
property of bilinearity and following rules for $i, \; j, \; k$:
\begin{eqnarray}
&&i^{2}=j^{2}=k^{2} =-1,\;\; jk=-kj = i, \nonumber\\
&&ki=-ik=j,\;\; ij=-ji= k.\nonumber
\end{eqnarray}
For $q=\alpha + \beta i+ \gamma j + \delta k$, the quaternionic conjugation
is defined as follows:
\begin{equation}
q^+=\alpha - \beta i -\gamma j - \delta k.
\end{equation}

Furthermore, we restrict ourselves to the set of quaternions
 ${\sf H}_{\mathbb R}\subset {\sf H}_{\mathcal C}$:
\begin{eqnarray*}
{\sf H}_{\mathbb R} =\{\zeta=\zeta^0 + {\rm i}(\zeta^1 i+ \zeta^2 j+
\zeta^3 k):\;{\rm i}^2=-1,\;\rm i\in{\mathbb C},\;
\zeta^0,\zeta^1,\zeta^2,\zeta^3  \in {\mathbb R} \},
\end{eqnarray*}
with the norm $\|\zeta\|^2$ given by
\begin{equation}
\|\zeta\|^2=\zeta {\zeta}^+=(\zeta^0)^2 - (\zeta^1)^2-(\zeta^3)^2 -
(\zeta^4)^2.
\end{equation}

Introducing a binary operation
\begin{equation}
\zeta\ast\eta=\big(\zeta+\eta\big)\big/\big(1 + \frac{K}{4}\zeta^+\eta \big),
\quad \zeta,\eta \in {\sf H}_{\mathbb R}, \label{QL}
\end{equation}
where $K$ is constant and $/$ denotes the right division, we find that the
set of quaternions ${\sf H}_{\mathbb R}$ with the binary operation $\ast$
forms a loop, which we denote by Q${\sf H}_{\mathbb R}$. The associator is found
to be
\begin{equation}
l_{(\zeta,\eta)}\xi= \big(1+\frac{K}{4}\zeta\eta^+ \big)\xi
\big/\big(1+\frac{K}{4}\zeta^+\eta \big).
\end{equation}

The loop Q${\sf H}_{\mathbb R}$ is related to the spacetime of constant
positive curvature $K>0$ ({\em de Sitter spacetime}), which is locally
characterized by the condition \cite{NS3}
\begin{equation}
R_{\mu\nu\lambda\sigma}= K(g_{\mu\lambda}g_{\nu\sigma} -
g_{\mu\sigma}g_{\nu\lambda}). \label{R}
\end{equation}
\end{example}

\subsection{Infinitesimal theory of smooth loops in brief}

It is well known that the infinitesimal theory of Lie groups arises from the
associativity of the operation
\[
a{\mathbf \cdot}(b{\mathbf \cdot} c)=(a{\mathbf \cdot} b){\mathbf
\cdot} c.
\]
In what follows, we will show that the quasiassociativity identities
\begin{equation}
R_c\circ R_b a=R_{R^a_c b}a, \quad L_a \circ L_b c = L_{(a{\mathbf
\cdot})}\circ l_{(a,b)}c
\end{equation}
lead to the infinitesimal theory of smooth loops. We start with the definition of a quasigroup of transformations introduced by Batalin \cite{Bat}.

\begin{definition} Let $\mathfrak M$ be an n-dimensional manifold and let the
continuous law of transformation be given by $x'=T_a x , \quad  x\in
{\mathfrak M}$, where \{$a^i$\} is the set of real parameters,
$i=1,2,\dots,r$. The set of transformations $\{T_a\}$ forms a
$r$-parametric quasigroup of
transformations (with right map on ${\mathfrak M}$), if \\
1) there exists a unit element that is common for all $x^i$ and corresponds to
$a^i=0:\;T_a x|_{a=0}=x$;\\
2) the modified composition law holds:
\[
T_aT_b x =T_{\varphi(b,a;x)} x;
\]
3) the left and right units coincide:
\[
\varphi(a,0;x)=a, \quad \varphi(0,b;x)=b;
\]
4) the modified law of associativity is satisfied:
\[
 \varphi(\varphi(a,b;x),c;x)=\varphi(a,\varphi(b,c;T_a x);x);
\]
4) the inverse transformation of $T_a$ exists: $x=T^{-1}_a x'$.
\end{definition}

The generators of infinitesimal transformations
\[
\Gamma_i=(\partial(T_ax)^i/\partial a^i)|_{a=0}
\partial/\partial x^\alpha \equiv R^\alpha_i\partial/\partial x^\alpha
\]
form a {\it quasialgebra} and obey the commutation relations
\begin{equation}
[\Gamma_i, \Gamma_j] =C^p_{ij}(x)\Gamma_p, \label{eq1}
\end{equation}
where $C^p_{ij}(x)$ are the structure functions satisfying the modified Jacobi
identity
\begin{eqnarray}
C^p_{ij,\alpha}R^\alpha_k + C^p_{jk,\alpha}R^\alpha_i
+C^p_{ki,\alpha}R^\alpha_j
+ C^l_{ij}C^p_{kl} + C^l_{jk}C^p_{il} + C^l_{ki}C^p_{jl}=0 . \label{eq2a}
\end{eqnarray}

{\begin{theorem}
Let the given functions $R^\alpha_i, \;C^p_{kj}$
obey the equations (\ref{eq1}), (\ref{eq2}). Then, locally the
quasigroup of transformations can be reconstructed as the solution of
the set of differential equations
\begin{eqnarray}
\frac{\partial\tilde x^\alpha}{\partial a^i} = R^\alpha_j(\tilde
x)\lambda^j_i(a;x),\quad
\tilde x^\alpha(0)=x^\alpha, \label{eq3} \\
\frac{\partial\lambda^i_j}{\partial a^p} -\frac{\partial\lambda^i_p}{\partial
a^j} + C^i_{mn}(\tilde x)\lambda^m_p\lambda^n_j=0, \quad
\lambda^i_j(0;x)=\delta^i_j \label{eq4}.
\end{eqnarray}
\end{theorem}}

Eq. (\ref{eq3}) is an analog of the Lie equation, and Eq.
(\ref{eq4}) is the generalized Maurer-Cartan equation.

The Batalin approach can be easily extended to the case of smooth loops
if we consider the action of the loop on itself. Indeed, for smooth loops, the modified associativity law can be written as follows:
\begin{equation} \varphi^i(\varphi (a,b),c)=\varphi^i(a,\tilde \varphi (b,c;a)), \quad
\varphi^i(a,0)=a^i, \quad \varphi^i(0,b)=b^i. \label{ID_1}
\end{equation}
The derivation of Eq.(\ref{ID_1}) in $e$ with respect to $c^j$  yields
\begin{equation}
\left.\frac{\partial\varphi^i(\varphi (a,b),c)}{\partial c^j} \right| _e =
\left.\frac{\partial\varphi^i(a,b)}{\partial b^k}
\frac{\partial\tilde\varphi^k(b,c;a)}{\partial c^j}\right| _e \label{Lie1}
\end{equation}
Let us introduce
\begin{equation}
\alpha^k_j(b;a)=\left.\frac{\partial \tilde\varphi^k(b,c;a)}{\partial c^j}
\right|_e,\quad \omega^i_p \alpha^p_j=\delta^i_j, \quad  q^i_p=
\left.\frac{\partial \varphi^i(a,b)}{\partial b^p} \right|_e. \label{Def}
\end{equation}
Substituting (\ref{Def}) in (\ref{Lie1}), we obtain an analog of Lie equation
(see Eq.(\ref{eq3}))
\begin{equation}
\frac{\partial\varphi^i(a,b)}{\partial b^j} =
q^i_p(\varphi(a,b))\omega^p_q(b;a). \label{Lie_1}
\end{equation}
Actually, the parametric operation $\varphi(b,c;a)$ is not independent, and it can be written in the terms of left translations as follows:
\[
\varphi(b,c;a) = L_b \circ l^{-1}_{(a,b)}c.
\]
This yields
\[
\omega^p_q(b;a) = [l_{(a,b),\ast,e}]^p_r \omega^r_q(b;e).
\]
Taking into account that $\omega^r_q(b;e)= \omega^r_q(b)$, we obtain the
generalized Lie equation deduced by Sabinin \cite{Sab1}
\begin{equation}
\frac{\partial\varphi^i(a,b)}{\partial b^j}  =
q^i_p(\varphi(a,b))[l_{(a,b),\ast,e}]^p_r \omega^r_q(b).
\end{equation}

\section{Principal Q-bundles}

\subsection{Basic definitions and examples}

Let $\mathfrak M$ be a manifold and $\langle Q, {\mathbf \cdot},
e\rangle$ be a smooth two-sided loop.  A {\em principal loop Q-bundle}
(a principal bundle with structure loop Q) is a triple $(P,\pi,
M)$, where $P$ is a manifolds and the following conditions hold:

(1) $Q$ acts freely on $P$ by the right map: $(p,a)\in P\times Q
\mapsto{} \tilde R_a p \equiv pa \in P$.

(2) $\mathfrak M$ is the quotient space of $P$ by the equivalence relation
induced by $Q$, $\mathfrak M=P/Q$, and the canonical projection  $\pi:
P\rightarrow M$ is a smooth map.

(3) $P$ is {\em locally trivial}, i.e. any $x\in  M$ has a
neighborhood $U$ and a diffeomorphism $\Phi$ of
$\pi^{-1}(U)\mapsto{} U\times Q$ such that, for any point $p\in
\pi^{-1}(U)$, it has the form $\Phi(p)=(\pi(p),\varphi(p))$, where
$\varphi$ is the map from $\pi^{-1}(U)$ to $Q$ satisfying
$\varphi(\tilde R_a p) = R_a \varphi(p)$.

We say that $P$ is a {\em total space} or {\em bundle space}, $Q$ is
the {\em structure loop}, $M$ is the {\em base} or {\em base
space} and $\pi$ is the {\em projection}. For any $x\in M$, the
inverse image $\pi^{-1}(x)$ is the {\em fiber} $F_x$ over $x$. Any
fiber is diffeomorphic to $Q$, and a transitive map of the loop $Q$
is defined on each fiber so that
\[
\tilde R_{b}\circ\tilde R_{a} p= \tilde R_{a{\mathbf \cdot}b}\tilde
r_{(a,b)}p,\quad
\]
where the map $\tilde r_{(a,b)}$ satisfies $\varphi(\tilde r_{(a,b)}p)=
r_{(a,b)}\varphi(p)$. The neutral element acts as the identity map.

Let us consider a covering $\{ U_\alpha\}_{\alpha \in J}$ of $ \mathfrak  M$
with $\bigcup_{\alpha \in J} U_\alpha =  \mathfrak M $, which can be chosen in
such a way that the restriction of the fibration to each open set $U_\alpha$ is
{\em trivializable }. This means that there exists a diffeomorphism
$\Phi_\alpha$: $\pi^{-1}(U_\alpha) \rightarrow{} U_\alpha\times Q$. The sets
$\{U_\alpha,\Phi_\alpha\}$ are called a family of {\em local trivializations}.
It is easy to prove the following \cite{N1,N2}.

{\begin{proposition} The right map on the fiber does not depend on
the trivialization: $\Phi^{-1}_\alpha
\big(\pi(p),\varphi_\alpha(\tilde R_q p) \big ) = \Phi^{-1}_\beta
\big(\pi(p),\varphi_\beta (\tilde R_q p) \big )$. \end{proposition}
}

{\begin{definition} Let $\{U_\alpha, \Phi_\alpha\}$ be a local
trivialization and $e\in Q$ a neutral element; then,
$\sigma_\alpha=\Phi^{-1}_\alpha(x,e),\; x\in U_\alpha$ is called the
{\em canonical local section} associated with this
trivialization.\end{definition}}

Let $u\in P$; then, one can describe the local section as follows:
$\sigma_\alpha = \tilde R^{-1}_{q_\alpha(u)}u$, where we set $q_\alpha(u)
= \varphi_\alpha (u).$

{\begin{proposition} The section $\sigma_\alpha$ does not depend on the choice of the point in the fiber.\end{proposition}}

{\begin{proof} Let $u,p\in P$ and $u = \tilde
R_{q_\alpha(u)}\sigma_\alpha$, $p = \tilde
R_{q_\alpha(p)}\sigma'_\alpha$. There exists an element $a\in Q$
such that $p=\tilde R_a u$. Then we have
\begin{eqnarray}
p= \tilde R_{q_{\alpha}(p)}\sigma'_{\alpha}= \tilde R_a u =\tilde
R_a \circ \tilde R_{q_\alpha(u)}\sigma_\alpha = &
\varphi^{-1}_\alpha\bigl(R_a \circ R_{q_\alpha(u)}
\varphi_\alpha(\sigma_\alpha) \bigr) \nonumber \\
= &\varphi^{-1}_\alpha\bigl(R_a\circ R_{q_\alpha(u)}e \bigr)= \tilde
R_{q_\alpha {\mathbf \cdot} a}\sigma_\alpha.
\end{eqnarray}
Taking into account that $q_{\alpha(p)} = q_{\alpha(u)}{\mathbf
\cdot}a$, we obtain
$R_{q_{\alpha}(p)}\sigma'_{\alpha}=R_{q_{\alpha}(p)}\sigma_{\alpha}$.
This yields $\sigma_\alpha=\sigma'_\alpha$, and, hence, the section
$\sigma_\alpha$ does not depend on the point of the fiber, while the
dependence on the base, $\sigma_\alpha=\sigma_\alpha(x)$, may still
hold. \end{proof}}

{\begin{proposition} Let the point $x$ lie in the intersection of the
neighborhoods $U_\alpha$ and $U_\beta$, $x\in U_\alpha\cap U_\beta$.
Then, the following formula holds:
\begin{eqnarray}
\sigma_\alpha(x)=\tilde R_{q_{\beta\alpha}}\sigma_\beta(x).
\label{sigma}
\end{eqnarray}\end{proposition}}

{\begin{proof} Let $u=\tilde R_{q_\alpha}\sigma_\alpha$, rewriting
this formula in the chart $U_\beta$, we obtain
\begin{eqnarray}
\Phi_\beta\bigl(\tilde R_{q_\alpha}\sigma_\alpha\bigr) =&\Phi_\beta\bigl(
\tilde R_{q_\beta}\sigma_\beta\bigr) =\bigl(\pi(u),R_{q_\beta}\varphi_\beta
(\sigma_\beta) \bigr)
\nonumber \\
=&\bigl(\pi(u),R_{q_{\beta\alpha}{\mathbf \cdot}q_\alpha}
\varphi_\beta(\sigma_\beta) \bigr)=  \bigl(\pi(u),R_{q_\alpha}\circ
R_{q_{\beta\alpha}}e \bigr) \nonumber \\
= &\bigl(\pi(u),R_{q_\alpha}\circ R_{q_{\beta\alpha}} \varphi_\beta(\sigma_\beta
)\bigr) =\Phi_\beta\bigl(\tilde R_{q_\alpha}\circ \tilde
R_{q_{\beta\alpha}}\sigma_\beta\bigr), \nonumber
\end{eqnarray}
where the relation $\varphi_\beta(\sigma_\beta)=e$ has been used.
Finally, we obtain $\sigma_\alpha(x)=\tilde
R_{q_{\beta\alpha}}\sigma_\beta(x)$ in the intersection
$U_\alpha\cap U_\beta$ . \end{proof}}

{\begin{definition} The family of maps $\{q_{\beta\alpha}(p) =
R^{-1}_{q_\alpha}q_\beta: \;\pi^{-1}( U_\alpha \cap
U_\beta)\rightarrow {} Q\}$, where $q_\alpha: =
{\varphi_\alpha(p)},\; q_\beta: = \varphi_\beta(p),\; p\in \pi^{-1}(
U_\alpha \cap U_\beta)$, is called the family of {\em transition
functions} of the bundle ${P({ M},Q)}$ corresponding to the
trivializing covering $\{U_\alpha,\Phi_\alpha\}_{\alpha \in J}$ of
$\mathfrak M$.\end{definition}}

{\begin{proposition} The transition functions $q_{\beta\alpha}(p)$
change under right translations as
\[
R_a q_{\beta\alpha}(p) = r_{(q_\alpha,a)}q_{\beta\alpha}(p)
\]
where $r_{(q_\alpha,a)}= R^{-1}_{q_\alpha{\mathbf \cdot}a}\circ R_a
\circ R_{q_\alpha}$ is the right associator.\end{proposition}  }

{\begin{proof} According to the definition, we have
\begin{eqnarray}
R_a q_{\beta\alpha}(p) = q_{\beta\alpha}(pa)
=R^{-1}_{q_\alpha{\mathbf \cdot} a} \circ R_a q_\beta
=R^{-1}_{q_\alpha{\mathbf \cdot} a} \circ R_a \circ
R_{q_\alpha}\circ R^{-1}_{q_\alpha}q_\beta=
r_{(q_\alpha,a)}q_{\beta\alpha}(p). \nonumber
\end{eqnarray}
\end{proof}}

We find that $q_{\beta\alpha}$ is changed under the right
translations, and, therefore, it depends on the points of the fiber:
$q_{\beta\alpha}=q_{\beta\alpha}(p),\; p\in F_x ,\;x \in U_\alpha
\cap U_\beta$. The transition functions obey the so-called {\it cocycle condition}:
\begin{eqnarray}
&&q_{\beta\alpha}(p){\mathbf \cdot}q_\alpha(p)= q_{\beta\gamma}(p)
{\mathbf \cdot} \bigl(q_{\gamma\alpha}(p) {\mathbf \cdot}
q_\alpha(p)\bigr)     \\
&&\forall x \in U_\alpha \cap U_\beta \cap U_\gamma \; {\rm and}\; \forall
p \in F_x , \nonumber
\end{eqnarray}
which may also be written as
\begin{eqnarray}
&&q_{\beta\alpha}{\mathbf \cdot}q_\alpha= q_{\beta\gamma} {\mathbf
\cdot} q_{\gamma\alpha}{\mathbf \cdot} \bigl(l_{(q_{\beta\gamma},
q_{\gamma\alpha})} q_\alpha \bigr). \nonumber
\end{eqnarray}

{\begin{definition} Let $P({  M},Q,\pi)$ be a principal loop
$Q$-bundle and $\{U_\alpha\}$ some covering of the base $\mathfrak
M$. Then we say that the system $\{ { \mathfrak
M},Q,\{U_\alpha\},q_{\alpha\beta} \}$ is a {\em principal coordinate
$Q$-bundle}\end{definition}}

Assume that there are two families of local sections, $\sigma_\alpha$ and
$\sigma'_\alpha$. These sections are linked by the right translations. This
implies that $\exists q_\alpha \in Q: \sigma'_\alpha =
R_{q_\alpha}\sigma_\alpha.$ Using the set $\{\sigma_\alpha\}$, let us
introduce the new transition functions $q'_{\alpha\beta}$. Then, the
following transformation law holds:
\[
q'_{\alpha\beta}=L^{-1}_{q_{\alpha}}\circ R_{q_\beta}q_{\alpha\beta}.
\]

\begin{example}
Principal $QS^2$ bundle over $S^1$. Here, we will
show that $S^3$ is the $QS^2$ bundle over $S^1$. We consider the
unit-three sphere $S^3$ as a subset of $\mathbb C^2$:
\begin{eqnarray}\label{eq7}
S^3 =\left\{(z_1,z_2)\in  \mathbb C^2:\; |z_1|^2+|z_2|^2=1 \right\}.
\end{eqnarray}
We define the map $\pi: S^3 \rightarrow S^1 $  by
\begin{eqnarray}\label{eq8}
\pi(z_1,z_2) =\bigg (\frac{\bar z_1 + z_1}{2\sqrt{1-|z_2|^2}}, \frac{ z_1 -
\bar z_1}{2i\sqrt{1-|z_2|^2}},0\bigg ).
\end{eqnarray}
Parameterizing $S^3$ by
\begin{eqnarray*}
z_1 = \cos\frac{\theta}{2}\,e^{i\psi_1}, \quad z_2
=\sin\frac{\theta}{2}\,e^{i\psi_2}
\end{eqnarray*}
where $0\leq \theta \leq \pi$ and $\psi_1,\psi_2 \in \mathbb R$, we get
\begin{eqnarray}\label{eq9}
\pi\big(\cos\frac{\theta}{2}\,e^{i\psi_1},
\sin\frac{\theta}{2}\,e^{i\psi_2}\big)= (\cos\psi_1,\sin\psi_1,0),
\end{eqnarray}
and one can see that $\pi$ indeed maps $S^3$ to $S^1$. By this means, one can see that locally
$S^3 \approx S^1 \times S^2$; however, this is not true globally since the spaces $S^1 \times S^2$ and $S^3$ have different second homotopy groups \cite{G3}:
\begin{eqnarray}\label{eq9a}
  \pi_2(S^1 \times S^2) =\mathbb Z, \quad  \pi_2(S^3) = 0.
\end{eqnarray}

Let $U_{-}=(-\pi,\pi)$ and $U_{+}=(0,2\pi)$ be an open covering of $S^1$. We define the
local trivializations, $ \Phi_{\pm}:  \pi^{-1}(U_{\pm})\rightarrow U_{\pm}\times QS^2 $,
as follows:
\begin{eqnarray}\label{eq10}
\Phi_{-}\big(\cos\frac{\theta}{2}\,e^{i\psi_1},
\sin\frac{\theta}{2}\,e^{i\psi_2}\big)
= \big(e^{i\psi_1}, \tan\frac{\theta }{2}\,e^{i(\psi_2-\psi_1)}\big) \\
\Phi_{+}\big(\cos\frac{\theta}{2}\,e^{i\psi_1},
\sin\frac{\theta}{2}\,e^{i\psi_2}\big) = \big(e^{i\psi_1},
\tan\frac{\theta +\pi/2}{2}\,e^{i(\psi_2-\psi_1)}\big).
\end{eqnarray}

Let us set
  \begin{eqnarray}\label{eq11}
  \zeta_{-} = \tan\frac{\theta }{2}\,e^{i(\psi_2-\psi_1)} ,\quad
  \zeta_{+} = \tan\frac{\theta +\pi/2}{2}\,e^{i(\psi_2-\psi_1)},
  \end{eqnarray}
then the transition function relative to our trivialization determined by
\begin{eqnarray}
   &q_{+-}: U_{-}\cap U_{+} \rightarrow QS^2
\end{eqnarray}
and 
\begin{eqnarray} \label{eq12}
      & \zeta_{+}= L_{q_{+-}} \zeta_{-} = \displaystyle\frac{q_{+-} + \zeta_{-}}{1 - \bar q_{+-}\zeta_{-}},
\end{eqnarray}
is given by $q_{+-} = e^{i(\psi_2-\psi_1)}$, and we have $(\cos\psi_1,\sin\psi_1,0)\rightarrow\, e^{i(\psi_2-\psi_1)}\in QS^2$.

Since the structure loop $QS^2$ acts on the standard fiber $S^2$ by left
translations (\ref{eq12}), we indeed have a $QS^2$ principal bundle over $S^1$.
For the right translatons of $QS^2$ on $S^3$, we find
\begin{eqnarray}\label{eq14}
    \tilde R_{\eta}(z_1,z_2)  = \Bigg(\frac{z_1|1-\bar \eta z_2/z_1|}{\sqrt{1+ |\eta|^2}},
    \frac{z_2(1+ \eta z_1/z_2)|1-\bar \eta z_2/z_1|}{(1-\bar \eta z_2/z_1)\sqrt{1+
    |\eta|^2}}\Bigg), \quad \eta\in QS^2.
\end{eqnarray}
\end{example}

\begin{example}
Principal $QS^2$ bundle over $S^2$. We construct
the principal $QS^2$ bundle over $S^2$ by taking
\begin{tabbing}
\qquad Base ${\mathfrak  M} = S^2$ \qquad \= with coordinates
$0 \leq \theta \leq \pi, \; 0 \leq \varphi < 2\pi$, \\ \kill \\
\qquad Fiber $QS^2$ \> with coordinates $ \zeta \in {\mathbb
C}$.
\end{tabbing}
We split $S^2$ into two hemispheres $H_{\pm}$ with
$H_{+}\cap H_{-}$ being  a thin strip parameterized by the
equatorial angle $\varphi$. Thus, locally the bundle can be described
as
\begin{tabbing}
\qquad $H_{-}\times QS^2$  \qquad \= with coordinates
$(\theta,\varphi,
\zeta_{-})$,\\ \kill \\
\qquad $H_{+}\times QS^2$ \qquad \= with coordinates
$(\theta,\varphi,
\zeta_{+})$. \\
\end{tabbing}
The transition functions must be elements of the loop $QS^2$ in order to give a
principal $Q$-bundle. We choose to relate the $H_{+}$ and $H_{-}$  fiber
coordinates in $H_{+}\cap H_{-}$  by the left map (\ref{eq12})
\begin{eqnarray}\label{eq15}
L_{q_{+-}} \zeta_{-} = \frac{q_{+-} + \zeta_{-}}{1 - \bar q_{+-}\zeta_{-}}.
\end{eqnarray}
Since $q_{+-}$ maps $S^2$ to $S^2$, it is classified by $\pi_2(S^2)= \mathbb Z$.
This implies that
\begin{eqnarray}
\zeta_{+} =\big( L_{q_{+-}}\big)^n\zeta_{-},
\label{eq16}
\end{eqnarray}
and the power $n$ of $L_{q_{+-}}$ must be an integer to give a well-defined
manifold. Let us define $q^n_{+-}$ by the following relation:
$\big( L_{q_{+-}}\big)^n = L_{q^n_{+-}}$; then, we get
\begin{eqnarray}
q^n_{+-} = \frac{\zeta_{+}(1- \zeta_{-}\overline\zeta_{+})
-\zeta_{-}(1- \overline\zeta_{-}\zeta_{+})}
{1 - |\zeta_{-}|^2|\zeta_{+}|^2}.
\end{eqnarray}
Setting $q_{+-}=e^{i\gamma}\tan(\theta/2)$ and using Eq.(\ref{eq16}) ,
we obtain $q^n_{+-}=e^{i\gamma}\tan(n\theta/2)$.
\end{example}

\subsubsection{Associated bundles}

In this section, we will discuss a method of constructing {\em quasigroup}
fiber bundles that are associated with some principle $Q$-bundle. Let $F$ be a
space on which $Q$ acts (on the left) as the loop of transformations. The basic
idea is to construct for a particular principal $Q$-bundle $P({\mathfrak M},Q)$
a fiber bundle with fiber $F$. First, we need the concept of a $Q$-{\em product}
of a pair of spaces on which $Q$ acts by the right mappings.

{\begin{definition} Let $X$ and $Y$ be any pair of right $Q$-spaces.
Then, the $Q$-{\em product} of $X$ and $Y$ is the space of orbits
of the $Q$-mappping on the Cartesian product $X\times Y. $
\end{definition}}

This definition implies that there is defined an equivalence relation on
$X\times Y$ in which $(x,y) \equiv (x',y')$ iff $\exists q\in Q$ such that $x'=
R_{\displaystyle q} x$ and $y' = R_q y$. The $Q$-product is denoted $X\times_Q
Y$.

{\begin{definition} Let $P({\mathfrak M},\pi,Q)$ be a principal
$Q$-bundle and let $F$ be a manifold on which $Q$ acts by left maps.
Define $P_F :=P\times_Q F$ and the map $Q$ on $P_F$ as follows
\[
(p,\xi)\rightarrow (R_a p, \tilde L^{-1}_a \xi), \quad p\in P,\; \xi \in F,
\]
where $\tilde L_a$ denotes the left map on $F$. Then, $E({\mathfrak
M},Q,F,P) =P\times_Q F/Q$ is the {\em associated} fiber bundle over
$\mathfrak M$ with fiber $F$.
\end{definition}
The differential structure over $E({\mathfrak M},Q,F,P)$ is defined
as follows. Let $\pi_E$ be the projection map $E\rightarrow
{\mathfrak M}: (u,\xi) \rightarrow \pi(u)$, that is induced by the
projection $\pi: P\times {\mathfrak M}\rightarrow \mathfrak M$. For
any $x\in {\mathfrak M}$, the set $\pi^{-1}_E (x)$ is called the {\em fiber}
in $E$ over $x$ and the $Q$-map on $\pi^{-1}_E (U)\times F$ is
defined by $(x,a,\xi)\rightarrow (x,ab,\tilde L^{-1}_b \xi)$, where
$(x,a,\xi)\in U\times Q\times F, \; b\in Q, U\subset {\mathfrak M}$
and $\pi^{-1}_E (U)$ is the open manifold in $E$: $\pi^{-1}_E (U)
\approx U\times F$. Then, the projection $\pi_E^{-1}$ is a differentiable
mapping $E\rightarrow {\mathfrak M}$. A map $\sigma: {\mathfrak
M}\rightarrow E$, such that $\pi_E\circ \sigma= {\rm Id}$ is the identity map ${\mathfrak M}\rightarrow {\mathfrak M}$, is called a
section of the fiber bundle $E({\mathfrak M},Q,F,P)$. }

\subsection{Connection, curvature and Bianchi identities}

\subsubsection{Connection on principal $Q$-bundles}

Let $P({  \mathfrak M},Q)$ be a principal $Q$-bundle over the manifold $
\mathfrak M$. For any $u\in P$, we denote a tangent space at $u$ by $T_u(P)$ (or
$T_u$), and the tangent to the fiber passing through $u$ by ${\mathcal V}_u$. We
call ${\mathcal V}_u$ a {\em vertical subspace}. It is generated by the right
translations on the fiber: $ u \rightarrow{}  \tilde R_a u, \; a\in Q, \; u\in
P$:
\begin{eqnarray}
X_u=\left.\frac{d}{dt}\tilde R_{a(t)} u\right |_{t=0}, \quad X_u\in
{\mathcal V}_u, \; a\in Q, \; u\in P. \label{ver0}
\end{eqnarray}

Let $\gamma(t)\in   \mathfrak M$ be a smooth curve. A {\em horizontal lift} of
$\gamma$ is a curve $\tilde\gamma(t)\in P$ such that
$\pi(\tilde\gamma(t))=\gamma(t)$. Evidently, to determine $\tilde\gamma(t)$,
it is sufficient to define at any point a tangent vector $\tilde X:
\pi_\ast \tilde X = X$, where $X$ is the tangent vector to $\gamma(t)$. A set
$\{ X\}$ is called a {\em horizontal subspace} ${\mathcal H}_u$.

{\begin{definition} A {\em connection} on a principal $Q$-bundle is
a smooth assignment to each point $p\in P$ of a horizontal subspace
${\mathcal H}_p$ of $T_p(P)$ such that:

(i) $ T_p={\mathcal V}_p \oplus{\mathcal H}_p$ (a direct sum).

(ii) The family of horizontal subspaces is invariant under the right
map of the loop $Q$; i.e., for any $a\in Q, \;p\in P$
\[
{\mathcal H}_{\tilde R_a p}= (\tilde R_a)_\ast {\mathcal H}_ p,
\]
where $\tilde R_{a\ast}$ denotes the right map by $a$ on a tangent
vector. \end{definition}}

The connection allows us to decompose any vector $Z\in T_p(P)$ in
the form $Z=X+Y$, where $X=hor Z \in {\mathcal H}_p$ is the
horizontal component of the vector $Z$ and $Y=ver Z \in {\mathcal
V}_p$ is the vertical one. The map $\varphi$ induces the map of the
vertical subspace ${\mathcal V}_p$ onto the tangent space to $Q$,
${\mathcal V}_p \stackrel{\varphi_\ast}\longrightarrow T_q(Q),\;
q=\varphi(p)$:
\begin{eqnarray}
{\tilde V}_p=\left.\frac{d}{dt}\tilde R_{a(t)} p\right |_{t=0}
\stackrel{\varphi_\ast}{\longrightarrow}
{\hat V}_q=\left.\frac{d}{dt}R_{ a(t)} q\right |_{t=0} .
\label{ver1}
\end{eqnarray}
Notice that the vector field $\hat V$ is the left quasi-invariant
vector field. Indeed, it can be written as $\hat V_q=(L_{q})_\ast
V_e$, where
\begin{eqnarray}
V_e=\left. \frac{da(t)}{dt}\right|_{t=0} \in T_e(Q).
\label{ver2}
\end{eqnarray}

\begin{definition} Let $ V_e \in T_e(Q)$. The vector field $\tilde
V$ connecting with ${V}_e$ by means of (\ref{ver1}), (\ref{ver2}) is
called a {\em fundamental} vector field.
\end{definition}

The notion of connection may be reformulated in the following way:

{\begin{definition} A connection form on a principal $Q$-bundle is a
vector-valued 1-form taking values at $T_e(Q)$, which satisfies:

(i) $\omega(X_p)= X_e$, where $X_p\in {\mathcal V}_p, \; X_e\in
T_e(Q)$ are determined according to {\rm (\ref{ver1}), (\ref{ver2})}.

(ii)$ \bigl(\tilde R^\ast_a \omega\bigr)X_p =
{\sf Ad}^{-1}_a(q) \omega(X_p)$, where $q=\varphi(p)$.

(iii) The horizontal subspace ${\mathcal H}_p$ is defined as a
kernel of $\omega$:
\[
{\mathcal H}_p=\{X_p \in T_p(P): \omega(X_p)=0 \}.
\]
\label{conn} \end{definition}}

A local 1-form taking values in $T_e(Q)$ can be
associated with the given connection form as follows. Let $\sigma: U\subset { \mathfrak M}\rightarrow{}
\sigma(U) \subset P,\; \pi\circ \sigma ={\rm id}$ be a local section of a
$Q$-bundle $Q\rightarrow{} P \rightarrow{} { \mathfrak  M}$ that is equipped
with a connection 1-form $\omega$. Define the {\em local
$\sigma$-representative} of $\omega$ to be the vector-valued 1-form $\omega^U$ (taking
values at $T_e(Q)$) on the open set $U\subset { \mathfrak  M}$ given
by $\omega^U:=\sigma^\ast \omega.$

{\begin{theorem} {\rm (On reconstruction of the connection form)}
For a given canonical {\sf Ad}-form $\tilde\omega$ defined on
$U\subset {  M}$ with values in $T_e(Q)$ and a given section $\sigma: U
\rightarrow{} \pi^{-1}(U)$, there exists one and only one connection
1-form $\omega$ on $\pi^{-1}(U)$ such that $\sigma^\ast
\omega=\tilde\omega$. \label{Rec}\end{theorem} }

{\begin{proof} Let $p_0=\sigma(x)$ and $Z\in T_{p_0}(P)$. We have
$Z=X_1 +X_2$, where $X_1:= (\sigma_\ast\circ\pi_\ast)Z$  and $X_2\in
{\mathcal V}_{p_0}, \; \pi_\ast X_2 = 0.$  Define at $p_0$ the
1-form $\omega$ to be the vector-valued 1-form given by
$\omega_{p_0} =\tilde\omega_x (\pi_\ast X) + \hat X_2$. The continuation
of the 1-form $\omega$ onto all points of the fiber is realized by
means of right translations, namely, $\forall p\in P,\;\exists a\in
Q: p=\tilde R_a p_0.$ This yields
\[
\omega_p\bigl((\tilde
R_a)_\ast X\bigr) = {\sf Ad}^{-1}_a(q_0)\omega_{p_0}(X), \quad
q_0=\varphi(p_0).
\]
It is easy to see that the 1-form so obtained satisfies all the conditions
of Definition \ref{conn}. \end{proof}}

This construction can be generalized to the whole base $\mathfrak M$ and leads
to the following theorem.

{\begin{theorem} Any smooth principal $Q$-bundle $(P,\pi,\mathfrak
M)$ has a connection. \end{theorem} }

{\begin{proof} Let $\tilde\omega$ be a differentiable 1-form on
$\mathfrak M$ with values in $T_e(Q)$ and $\{ U_\alpha,
\Phi_\alpha\}$ a family of local trivializations associated with the
covering $\{ U_\alpha\}$ of $\mathfrak M$. For each local
trivialization, there exists the local connection form
$\omega_\alpha$ on $\pi^{-1}(U_\alpha)$ obtained from $\tilde\omega$
by a choice of a section $\sigma_\alpha$ over $U_\alpha$ (see
Theorem \ref{Rec}). Let $\lambda_\alpha$ be a partition of unity
subordinate to the covering $\{U_\alpha\}$, $\lambda = 0$ outside of
$U_\alpha$ and $\sum_\alpha \lambda_\alpha (x) =1$ at an arbitrary
point  $x \in  \mathfrak M$. We form
\[
\omega = \sum_\alpha (\lambda_\alpha \circ \pi) \omega_\alpha.
\]
This is a 1-form on $P$ with values in $T_e(Q)$ and for each $\alpha$
we have $\omega\left|_{{}_{U_\alpha}}\right. =
(\lambda_\alpha\circ\pi) \omega_\alpha$. It is easily seen that
$\omega$ satisfies the conditions of the Definition \ref{conn} and
thus this is a connection form on $P$. \end{proof}}

\subsubsection{Connection in the local trivialization }

We consider a principal Q-bundle $(P,\pi, M)$ over $\mathfrak M$ and fix a
trivializating covering $\{U_\alpha,\Phi_\alpha\}_{\alpha \in J}$. Let
$\overline{\rm Id}:  U_\alpha\rightarrow{} U_\alpha\times Q$ by $x
\rightarrow{} (x,e)$. A trivialization $\Phi_\alpha$ defines a {\em canonical}
section $\sigma_\alpha$ by the equation

\[
\sigma_\alpha =
\Phi^{-1}_\alpha \circ\overline{\rm Id},
\]
and vice versa. Further, we denote by $\omega_{C}$ a canonical $\sf Ad$-form.

{\begin{definition} Let $\omega_\alpha = \sigma_\alpha^\ast \omega$,
where $\omega$ is the connection form. The form $\omega_\alpha$ on
$U_\alpha$ is called the {\em connection form in the local
trivialization} $\{U_\alpha,\Phi_\alpha\}.$ \end{definition}}

{\begin{theorem} Let $\{U_\alpha, \Phi_\alpha\}_{\alpha \in J}$ be a
family of local trivializations for $P$ with $\bigcup_{\alpha\in J}
U_\alpha = M$; then, on $U_\alpha\cap U_\beta$, the local connection forms
$\omega_\alpha$ and $\omega_\beta$, corresponding to the same
connection $\omega$ on $P$, are related by
\begin{eqnarray}
\omega_\beta = {\sf Ad}^{-1}_{q_{\alpha\beta}}(q_{\beta\alpha})
\omega_\alpha + l_{(q_{\beta\alpha},q_{\alpha\beta})\mathbf \ast}
\theta_{\alpha\beta}, \label{local}
\end{eqnarray}
where the $q_{\alpha\beta}$ are the transition functions, and
$\theta_{\alpha\beta}=q^{\ast}_{\alpha\beta}\omega_{C}$ denotes the
pullback on $U_\alpha\cap U_\beta$ of the canonical 1-form
$\omega_{C}$ on $Q$. Vice versa, for any set of the local forms $\{
\omega_\alpha\}_{\alpha \in J}$ satisfying {\rm(\ref{local})}, there
exists the unique connection form $\omega$ on $P$ generating this
family of local forms, namely
$\omega_\alpha=\sigma^\ast_\alpha\omega,\; \forall\alpha \in J$.
\end{theorem}}

{\begin{proof}
1. The direct theorem. Let $x\in U_\alpha\cap
U_\beta$. Applying (\ref{sigma}), we obtain $\sigma_\beta(x) =
\tilde R_{q_{\alpha\beta}}\sigma_\alpha(x),\; \forall x\in
U_\alpha\cap U_\beta$. The map $(\sigma_\beta(x))_\ast$ transforms
any vector $X\in T_x(U_\alpha\cap U_\beta)$ into
$(\sigma_\beta)_\ast X \in T_{\sigma_\beta(x)}(P)$. Using the
Leibniz formula (see \cite{KN1}, Ch. I), we get
\[
(\sigma_\beta)_\ast = (\tilde  R_{q_{\alpha\beta}})_\ast
(\sigma_\alpha)_\ast + (L_{q_{\beta\alpha}})_\ast (q_{\alpha\beta})_\ast
\]
where $q_{\beta\alpha}:=\varphi_\beta(\sigma_\alpha)$. Applying $\omega$ to
both sides of this relation, we find
\begin{eqnarray}
\omega_\beta(X):=\omega((\sigma_\beta)_\ast X)&& =
\omega\bigl((\tilde R_{q_{\alpha\beta}})_\ast (\sigma_\alpha)_\ast
X\bigr)
+\omega\bigl((L_{q_{\beta\alpha}})_\ast (q_{\alpha\beta})_\ast X\bigr) \nonumber\\
&&={\sf
Ad}^{-1}_{q_{\alpha\beta}}(q_{\beta\alpha})\,\omega_{\alpha}(X) +
l_{(q_{\beta\alpha},q_{\alpha\beta})\ast}
(q^\ast_{\alpha\beta}\,\omega_{C})(X)\nonumber \\
&&={\sf
Ad}^{-1}_{q_{\alpha\beta}}(q_{\beta\alpha})\,\omega_{\alpha}(X) +
l_{(q_{\beta\alpha},q_{\alpha\beta})\ast}
\theta_{\alpha\beta}\,\omega_{C}(X).
\end{eqnarray}

2. The inverse theorem. Let us define the 1-form $\tilde\omega$ as
follows:
\begin{eqnarray}
\tilde\omega={\sf Ad}^{-1}_{q_\alpha}(e)(\pi^\ast \omega_\alpha)
+ q^\ast_\alpha\, \omega_{C}, \quad q_\alpha:= \varphi_\alpha(p).
\label{L1}
\end{eqnarray}
Let $X\in T_u(P)$ be an arbitrary vector and $u=\sigma_\alpha(\pi(x))$.
Decompose $X$ into horizontal $Y$ and vertical $Z$ components:
\[
X=Y+Z, \quad Y=(\sigma_\alpha)_\ast (\pi_\ast X), \quad \pi_\ast Z =0.
\]
This yields
\begin{eqnarray}
\tilde\omega(X)&&= {\sf Ad}^{-1}_{q_\alpha}(e)\omega_\alpha(\pi_\ast
X ) + (q^\ast_\alpha\, \omega_{C})(X) ={\sf
Ad}^{-1}_{q_\alpha}(e)\omega((\sigma_\alpha)_\ast \pi_\ast X )
+\omega_{C} ((q_\alpha)_\ast X)\nonumber\\
&&= (\tilde R^\ast_{q_\alpha} \sigma^\ast_\alpha) \omega(\pi_\ast Y)
+ \omega_{C}((q_\alpha)_\ast Z) =
\omega((\sigma_\alpha)_\ast \pi_\ast Y )  + \omega_{C} ((q_\alpha)_\ast Z)\nonumber\\
&&=\omega(Y) + \hat Z = \omega(Y)+ \omega(Z) = \omega(X),
\end{eqnarray}
and one sees that $\tilde\omega = \omega$ at any point of the
section $\sigma_\alpha$. So, these forms are transformed in the same
way under the right translations and therefore coincide on
$\pi^{-1}(U_\alpha)$. \end{proof}}

{\begin{corollary} For arbitrary sections $\sigma_1$ and $\sigma_2$
such that $\sigma_2 = R_q \sigma_1$ and $\omega_1 = \sigma^\ast_1
\omega, \;\omega_2 = \sigma^\ast_2 \omega$, the following relation
holds:
\[
\omega_2 = {\sf Ad}^{-1}_q(q_1)\omega_1 + l_{(q_1,q))\ast} (q^\ast
\omega_{C}),
\]
where $q_1:=\varphi(\sigma_1)$. \end{corollary}}

\subsection{Covariant derivative. Curvature form}

Let $\{x^\mu, y^i\}$ be a local coordinate system in the neighborhood
$\pi^{-1}(U_\alpha)$, where the $x^\mu$ are coordinates in $U_\alpha \in {M}$ and
the $y^i$ are coordinates in the fiber. Locally $\pi^{-1}(U_\alpha)$ can be
presented as a direct product $U_\alpha\times Q$. The connection form
can be written in the form: $\omega=\omega^i L_i$, with $\{L_i\}$ being the
basis of left fundamental fields and $\{\omega^i\}$ the basis of 1-forms.
Taking into account Eq. (\ref{L1}), we find that in the coordinates
$\{x^\mu, y^i\}$,
\begin{eqnarray}
\omega^i= \bigl({\sf Ad}^{-1}_{\displaystyle y}(e)\bigr)^i_j
A^j_\mu(x)dx^\mu + \omega^i{}_j(y) dy^j
\label{om}
\end{eqnarray}
where $A^i_\mu(x)dx^\mu = \pi^{\ast}(\omega^i) $ and $\omega^i{}_j(y) dy^j
= \bigl(L^{-1}_{*}\bigr)^i_j dy^j$.

{\begin{definition} A covariant derivative $D_\mu$ in the principal
$Q$-bundle is defined as follows:
\[
D_\mu = \partial_\mu - A^i_\mu (x)\bar L_i,
\]
where the $\bar L_i=(R_\ast)^j_i \partial/\partial y^j$ are generators of the
left translations (right quasi-invariant vector
fields).\end{definition} }

Now, we show that $\omega(D_\mu) = 0$. Indeed, Eq. (\ref{om}) implies
\[
\omega(D_\mu)= [\bigl({\sf Ad}^{-1}_{\displaystyle y}(e)\bigr)^i_j
A^j_\mu(x) - \omega^i{}_p A^j_\mu (R_\ast)^p_j]L_i.
\]
Noting that $\omega^i{}_p= (L^{-1}_\ast)^i_ p$, we obtain
$\omega^i{}_p(R_\ast)^p_j = ({\sf Ad}^{-1}_q (e))^i_j$, and hence
$\omega(D_\mu) =0$.

The computation of the commutator $[D_\mu,D_\nu]$ yields
\[
[D_\mu,D_\nu]=(\partial_\nu A^i_\mu - \partial_\mu A^i_\nu)\bar L_i
+ A^i_\mu A^j_\nu [\bar L_i, \bar L_j],
\]
and introducing $[\bar L_i, \bar L_j] = C^p_{ij}(y)\bar L_p$, we get
\begin{eqnarray}
[D_\mu,D_\nu]=- F^i_{\mu\nu}\bar L_i,\\
F^i_{\mu\nu}: = \partial_\mu A^i_\nu - \partial_\nu A^i_\mu
- A^j_\mu A^p_\nu C^i_{jp}.
\end{eqnarray}

{\begin{definition} Let $\Psi$ be a vector-valued $r$-form in the
principal $Q$-bundle. A $(r+1)$-form $D\Psi$ defined by
\begin{eqnarray}
D\Psi(X_1,X_2,\dots,X_{r+1})=d\Psi(hor X_1,\dots,hor X_{r+1})
\end{eqnarray}
is called a {\em covariant differential} of the form $\Psi$.
\end{definition}}

{\begin{definition} A vector-valued 2-form $\Omega(X,Y)$ defined as
\begin{eqnarray}
\Omega(X,Y) = D\omega (X,Y) =d\omega (hor X, hor Y)
\end{eqnarray}
where $\omega$ is a connection form, is called a {\em curvature form}.
\end{definition}}

{\begin{lemma} Let $X,Y$ be horizontal fields, then the following
relation holds:
\[
\omega([X,Y]) = -2\Omega(X,Y).
\]
\label{lem2}\end{lemma} }

{\begin{proof} Applying the exterior differentiation to the 1-form
$\omega$, we obtain
\[
2d\omega(X,Y)=X\omega(Y) - Y\omega(X) -\omega([X,Y]).
\]
As $X,Y \in {\mathcal H}_p$, then $\omega(X)=\omega(Y)=0$. This
yields $\omega([X,Y]) = -2\Omega(X,Y)$. \end{proof}}

{\begin{corollary} The curvature form can be defined as follows:
\[
\Omega(X,Y)=-\frac{1}{2}\omega([hor\bar Y,hor\bar Y])
\]
where $\bar X, \bar Y$ are any continuations of  the vectors $X,Y
\in T_p(P)$, respectively.\end{corollary} }

\begin{corollary}
The 2-form $\Omega$ is $\sf Ad$-form and is
transformed under the right translations as
\begin{eqnarray}
(\tilde R^\ast_a \Omega)(X,Y)= {\sf Ad}^{-1}_a(q)\,\Omega \vert_p
(X,Y), \; {\rm where} \; q=\varphi(p).
\end{eqnarray}
\end{corollary}

{\begin{theorem} The curvature form $\Omega$ satisfies the structure
equation
\begin{eqnarray}
\Omega = d \omega + \frac{1}{2}{[\omega,\omega]}.
\label{str}
\end{eqnarray}
\end{theorem}}

{\begin{proof} We are going to prove (\ref{str}), considering all
possible pairs $X,Y$.

1. Let $X,Y$ be any two horizontal vectors then $\omega(X)=\omega(Y)=0$ and
we find
\[
d\omega(X,Y) + \frac{1}{2}{[\omega(X),\omega(Y)]} = d\omega(hor X, hor Y) =
\Omega(X,Y)
\]

2. Let $X,Y$ be any two vertical vectors. Without loss of generality, one can
assume that $X,Y$ are fundamental vector fields. This implies $\omega(X) =
\hat X,\; \omega (Y) =\hat Y, \; \hat X, \hat Y \in T_e(Q)$.
The computation yields
\begin{eqnarray}
2 d \omega(X,Y)= X\omega(Y) - Y\omega(X) -\omega([X,Y])
= X(\hat Y) - Y(\hat X) -\omega([X,Y])   \\
= -[\hat X, \hat Y] =- [\omega(X),\omega(Y)] \Longrightarrow
d \omega(X,Y) + \frac{1}{2}[\omega(X),\omega(Y)] =0.
\end{eqnarray}
As $X,Y$ are the vertical fields, then $\Omega(X,Y)= d\omega(hor X,
hor Y) =0$.

3. Let $X$ be a horizontal vector field and $Y$ be a vertical one (or
fundamental vector field). From the one side, we have
\[
\Omega(X,Y)=d\omega(hor X, hor Y) =0,
\]
and from the other side,
\[
d\omega(X,Y)=X\omega(Y) - Y\omega(x) -\omega([X,Y]) =0,
\]
since the commutator $[X,Y]$ is the horizontal vector field.
\end{proof}}

The last statement arises from the following lemma.

{\begin{lemma}  Let $Z$ be a fundamental vector field and $X$ be a
horizontal vector field; then, the commutator $[X,Z]$ is a horizontal
vector field. \end{lemma}  }

{\begin{proof} The fundamental vector field is induced by the right
translations $\tilde R_{a(t)}$. The commutator can be defined (see,
e.g., \cite{KN1}) by
\[
[X,Z] = \lim_{t\rightarrow{} 0}\frac{1}{t}((\tilde R_{a(t)})_\ast X - X).
\]
If $X$ is a horizontal vector field, then $(\tilde R_{a(t)})_\ast X$
is a horizontal vector field, and, hence, $[X,Z]$ is also
horizontal. \end{proof}}

Let us compute $\Omega_{\mu\nu} = 2\Omega(D_\mu,D_\nu)$. The result is given by
\begin{eqnarray}
\Omega_{\mu\nu}= -\omega([D_\mu,D_\nu])
= F^i_{\mu\nu}\omega(\bar L_i)= F^i_{\mu\nu}
\bigl({\sf Ad}^{-1}_q(e)\bigr)^j_i \hat L_j \nonumber
\label{Om}
\end{eqnarray}
where $\{\hat L_j =\omega(L_j)\}$ is the basis of the left quasi-invariant
vector fields at $T_e(Q),\; q=\varphi(p),\; p\in P$. Choosing the family of
 local sections $\sigma_\alpha$ associated with the trivialization
$U_\alpha, \Phi_\alpha$ and taking into account that
$\varphi_\alpha(\sigma_\alpha)=e$ where $\varphi_\alpha$ is the restriction of
$\Phi_\alpha$ on $\pi^{-1}(U_\alpha)$, we obtain
\begin{eqnarray}
\pi^\ast \Omega\vert_\alpha =\frac{1}{2}F_{\mu\nu}dx^\mu \wedge dx^\nu.
\end{eqnarray}
introducing $F_{\mu\nu}= F^i_{\mu\nu}L_Œ$.

Since $\Omega$ is a $\sf Ad$-form, the following transformation law is true:
\begin{eqnarray}
\Omega\vert_\beta = {\sf Ad}^{-1}_{q_{\alpha\beta }}
(q_{\beta\alpha})\Omega\vert_\alpha.
\end{eqnarray}

{\begin{theorem} Bianchi identity: $D\Omega =0$. \end{theorem}}

{\begin{proof} It is sufficient to show that $d\Omega(X,Y,Z) =0$, if
$X,Y,Z$ are horizontal vector fields. Applying the exterior
derivative  to (\ref{str}), we obtain $d\Omega(X,Y,Z) =0$, if
$X,Y,Z$ are horizontal vector fields. \end{proof}}

\section*{Acknowledgements}

This work was partly supported by SEP-PROMEP (Grant No.
103.5/04/1911).

\end{document}